\documentstyle[amscd, leqno]{amsart}
\theoremstyle{plain}
\newtheorem{thm}{Theorem}[section]
\newtheorem{Mthm}[thm]{Main Theorem}
\newtheorem{prop}[thm]{Proposition}
\newtheorem{lem}[thm]{Lemma}

\newtheorem{cor}[thm]{Corollary}
\theoremstyle{definition}
\newtheorem{defn}[thm]{Definition}
\newtheorem{exmp}[thm]{Example}
\theoremstyle{remark}
\newtheorem{rem}[thm]{Remark}
\makeatletter

\@addtoreset{equation}{section}
\makeatother

\begin{document}
\title[Real $K3$ surfaces without real points]
{Real $K3$ surfaces without real points,
\\
equivariant determinant of the Laplacian,
\\
and the Borcherds $\Phi$-function }
\author{Ken-Ichi Yoshikawa}
\address{
Graduate School of Mathematical Sciences,
University of Tokyo,
3-8-1 Komaba, Tokyo 153-8914, JAPAN}
\email{yosikawa@@ms.u-tokyo.ac.jp}
\thanks{The author is partially supported 
by the Grants-in-Aid for Scientific Research for 
Encouragement of Young Scientists (B) 16740030, JSPS}

\begin{abstract}
We consider an equivariant analogue of 
a conjecture of Borcherds. Let $(Y,\sigma)$ be
a real $K3$ surface without real points.
We shall prove that the equivariant determinant 
of the Laplacian of $(Y,\sigma)$ with respect to 
a $\sigma$-invariant Ricci-flat K\"ahler metric 
is expressed as the norm of the Borcherds 
$\Phi$-function at the ``period point''.
Here the period of $(Y,\sigma)$ is not the one
in algebraic geometry. 
\end{abstract}

\maketitle


\section
{\bf Introduction }
\par
Let $Y$ be an algebraic $K3$ surface defined 
over the real number field $\Bbb R$. 
Let $\sigma\colon Y\to Y$ be the anti-holomorphic 
involution on $Y$ induced by 
the complex conjugation. 
Denote by 
${\Bbb Z}_{2}=\langle\sigma\rangle$ the group of 
order $2$ of $C^{\infty}$ diffeomorphisms of $Y$
generated by $\sigma$. 
Recall that a point of $Y$ is real 
if it is fixed by $\sigma$.
\par
By \cite{Ya}, there exists a $\sigma$-invariant 
Ricci-flat K\"ahler metric $g$ on $Y$
with K\"ahler form $\omega_{g}$. 
Since $Y$ is defined over $\Bbb R$, 
there exists a nowhere vanishing holomorphic 
$2$-form $\eta_{g}$ on $Y$ such that 
$$
\eta_{g}\wedge\overline{\eta}_{g}=2\omega_{g}^{2},
\qquad
\sigma^{*}\eta_{g}=\overline{\eta}_{g}.
$$
Notice that the choice of $\eta_{g}$ is unique
up to a sign.
We identify $\omega_{g}$ and $\eta_{g}$ with
their cohomology classes.
\par
Let ${\Bbb L}_{K3}$ be the $K3$ lattice, which is
an even unimodular lattice with signature $(3,19)$. 
Then $H^{2}(Y,\Bbb Z)$ equipped with the cup-product 
is isometric to ${\Bbb L}_{K3}$.
By \cite{Ni1} or \cite{DIK00},
there exists an isometry of lattices
$\alpha\colon H^{2}(Y,{\Bbb Z})\cong{\Bbb L}_{K3}$ 
such that the point 
$[\alpha(\omega_{g}+\sqrt{-1}{\rm Im}\,\eta_{g})]
\in{\Bbb P}({\Bbb L}_{K3}\otimes{\Bbb C})$ lies in
the period domain for Enriques surfaces.
\par
Let $\Delta_{Y,g}$ be the Laplacian of $(Y,g)$ 
acting on $C^{\infty}(Y)$. 
Following \cite{Bi} and \cite{KR},
one can define the equivariant determinant 
of the Laplacian $\Delta_{Y,g}$ with respect to 
the anti-holomorphic ${\Bbb Z}_{2}$-action on $Y$. 
Notice that $\sigma$ acts on the vector space 
$C^{\infty}(Y)$ while it does not 
act on the vector space of $C^{\infty}$ $(p,q)$-forms 
on $Y$ unless $p=q$. 
Denote by 
${\rm det}_{{\Bbb Z}_{2}}^{*}\Delta_{Y,g}(\sigma)$ 
the equivariant determinant of 
the Laplacian $\Delta_{Y,g}$ 
with respect to $\sigma$. (See Sect.\,4.2.)
\par
Recall that Borcherds \cite{Bo1}
constructed a very interesting
automorphic form on the period domain 
for Enriques surfaces,
which is called the Borcherds $\Phi$-function and 
is denoted by $\Phi$. 
Let $\|\Phi\|$ denote the Petersson norm of $\Phi$. 
Then $\|\Phi\|^{2}$ is a $C^{\infty}$ function 
on the period domain for Enriques surfaces,
which is invariant under the complex conjugation of
the period domain. Our result is the following:

\begin{Mthm} 
There exists an absolute constant $C>0$ such that
for every real $K3$ surface without real points
$(Y,\sigma)$ and for every $\sigma$-invariant 
Ricci-flat K\"ahler metric $g$ on $Y$ 
with volume $1$,
$$
{\rm det}_{{\Bbb Z}_{2}}^{*}\Delta_{Y,g}(\sigma)
=C\,\|\Phi([\alpha(\omega_{g}+
\sqrt{-1}{\rm Im}\,\eta_{g})])\|^{\frac{1}{4}}.
$$
\end{Mthm}

Notice that the point
$[\alpha(\omega_{g}+\sqrt{-1}{\rm Im}\,\eta_{g})]$ 
is {\it not} the period of the marked $K3$ surface 
$(Y,\alpha)$,
because $\omega_{g}+\sqrt{-1}{\rm Im}\,\eta_{g}$ 
is not a holomorphic $2$-form on $Y$.
Since $\omega_{g}$ is the K\"ahler form of $(Y,g)$,
the Main Theorem 1.1 may be regarded as 
a symplectic analogue of \cite[Th.\,8.3]{Yo}.
A typical example of a real $K3$ surface 
without real points is the quartic surface 
of ${\Bbb P}^{3}({\Bbb C})$ defined by
the equation 
$z_{0}^{4}+z_{1}^{4}+z_{2}^{4}+z_{3}^{4}=0$.
\par
To prove the Main Theorem 1.1, we consider 
an equivariant analogue of the conjecture of
Borcherds:
Let $X$ be the differentiable manifold underlying
a $K3$ surface.
In \cite[Example 15.1]{Bo2}, Borcherds conjectured 
that the regularized determinant of the Laplacian, 
regarded as a function on the moduli space 
of Ricci-flat metrics on $X$ with volume $1$, 
coincides with the automorphic form 
on the Grassmann $G({\Bbb L}_{K3})$ associated to 
the elliptic modular form
$E_{4}(\tau)/\Delta(\tau)$; 
it is worth remarking that
the regularized determinant of the Laplacian 
of a Ricci-flat $K3$ surface can be regarded 
as an analytic torsion of certain elliptic complex 
\cite{KW}.
\par
As an equivariant analogue 
of the Borcherds conjecture,
we shall compare the following two functions 
on the space of $\sigma$-invariant Ricci-flat 
metrics on $X$;
one is the equivariant determinant of the Laplacian, 
and the other is the pull-back of the norm 
of the Borcherds $\Phi$-function via 
the ``period map''. 
(See Sect.\,3.4 for the definition
of the period map.)
It is a trick of Donaldson
\cite{DIK00}, \cite{D2} that
relates these two objects:
Let $(I,J,K)$ be a hyper-K\"ahler structure 
on $(X,g)$ with $Y=(X,J)$. Then $\sigma$ is 
holomorphic with respect to another complex 
structure $I$, while $\sigma$ is anti-holomorphic
with respect to the initial complex structure $J$.
We shall show that the equivariant determinant 
of the Laplacain of $(Y,\sigma)$ coincides with 
the equivariant analytic torsion of $(X,I,\sigma)$. 
(See Sect.\,3.3 and Sect.\,4.) 
After this observation, 
the Main Theorem 1.1 is a consequence of our result
\cite[Main Theorem and Th.\,8.2]{Yo}.
\par
This note is organized as follows. 
In Sect.\,2, we recall the notion of hyper-K\"ahler 
structure on a $K3$ surface.
In Sect.\,3, we recall the trick of Donaldson.
In Sect.\,4, we study equivariant determinant of the 
Laplacian as a function on the space of 
$\sigma$-invariant Ricci-flat metrics 
on a $K3$ surface, and we prove the Main Theorem.
\par
We thank Professors Jean-Michel Bismut and Sachiko Saito
for helpful discussions on the subject of this note. 
This note is inspired by \cite{Bost}.


\section
{\bf $K3$ surfaces and hyper-K\"ahler structures}
\par

\subsection
{}
{\bf $K3$ surfaces}
\par
A compact, connected, smooth complex surface is 
a {\it $K3$ surface} if it is simply connected 
and has trivial canonical line bundle.
Every $K3$ surface is diffeomorphic to 
a smooth quartic surface 
in ${\Bbb P}^{3}({\Bbb C})$ 
(cf. \cite[Chap.\,8 Cor.\,8.6]{BPV}).
Throughout this note, $X$ denotes the
$C^{\infty}$ {\it differentiable manifold underlying 
a $K3$ surface},
and $X$ is equipped with the orientation 
as a complex submanifold of ${\Bbb P}^{3}({\Bbb C})$.
For a complex structure $I$ on $X$, 
$X_{I}$ denotes the $K3$ surface $(X,I)$. 
\par
Let $\Bbb U$ be the lattice of rank $2$ associated 
with the symmetric matrix $\binom{0\,1}{1\,0}$,
and let ${\Bbb E}_{8}$ be the root lattice 
of the simple Lie algebra of type $E_{8}$. 
We assume that ${\Bbb E}_{8}$
is {\it negative-definite}. 
The even unimodular lattice with signature $(3,19)$
$$
{\Bbb L}_{K3}
:=
{\Bbb U}\oplus{\Bbb U}\oplus{\Bbb U}\oplus
{\Bbb E}_{8}\oplus{\Bbb E}_{8}
$$
is called {\it the $K3$ lattice}. 
Then $H^{2}(X,{\Bbb Z})$ equipped with 
the cup-product $\langle\cdot,\cdot\rangle$,
is isometric to ${\Bbb L}_{K3}$ 
(cf. \cite[Chap.\,8, Prop.\,3.2]{BPV}).

\subsection
{}
{\bf Hyper-K\"ahler structures on $X$}
\par
In this subsection, we recall Hitchin's result 
\cite{H}.
Let $\cal E$ be the set of all Ricci-flat metrics 
on $X$ with volume $1$. 
For every complex structure $I$ on $X$, 
there exists a K\"ahler metric on $X_{I}$ 
by \cite[Chap.\,8, Th.\,14.5]{BPV}. 
For every K\"ahler class $\kappa$ on $X_{I}$, 
there exists by \cite{Ya}
a unique Ricci-flat K\"ahler form on $X_{I}$ 
representing $\kappa$. 
Hence 
${\cal E}\not=\emptyset$. For $g\in{\cal E}$, 
let $dV_{g}$ denote the volume element of $(X,g)$. 
Then $\int_{X}dV_{g}=1$ by our assumption.

\begin{defn}
A complex structure $I$ on $X$ is {\it compatible}
with $g\in{\cal E}$ if $g$ is a K\"ahler metric 
on $X_{I}$, i.e., $I$ is parallel with respect 
to the Levi-Civita connection of $(X,g)$.
For $g\in{\cal E}$, let ${\cal C}_{g}$ denote 
the set of all complex structures on $X$ 
compatible with $g$.
\end{defn}

Let $g\in{\cal E}$.
By Hitchin 
\cite[Sect.\,2, (i) $\Leftrightarrow$ (iii)]{H}, 
we get ${\cal C}_{g}\not=\emptyset$.
For $I\in{\cal C}_{g}$, we define 
a real closed $2$-form $\gamma_{I}$ on $X$ by
\begin{equation}
\gamma_{I}(u,v):=g(Iu,v),
\qquad
u,v\in TX.
\end{equation}
Then $\gamma_{I}$ is a Ricci-flat K\"ahler form 
on $X_{I}$ such that 
$$
\gamma_{I}^{2}=2dV_{g}.
$$

\begin{defn}
Let $I,J,K\in{\cal C}_{g}$. 
The ordered triplet $(I,J,K)$ is called 
a {\it hyper-K\"ahler structure} on $(X,g)$ if 
\begin{equation}
IJ=-JI=K.
\end{equation}
\end{defn}

Let 
$*_{g}\colon\bigwedge^{p}T^{*}X\to
\bigwedge^{4-p}T^{*}X$ 
be the Hodge star-operator on $(X,g)$. 
Since $\dim_{\Bbb R}X=4$, 
we have $*_{g}^{2}=1$ on $\bigwedge^{2}T^{*}X$. 
Recall that a $2$-form $f$ on $X$ is 
{\it seld-dual} with respect to $g$ if $*_{g}f=f$.
Let ${\cal H}^{2}_{+}(g)$ be the real vector space 
of self-dual, real harmonic $2$-forms on $(X,g)$. 
Every vector of ${\cal H}^{2}_{+}(g)$ is 
parallel with respect to the Levi-Civita 
connection by \cite{H}.

\begin{thm}
Let $I\in{\cal C}_{g}$, and let $\eta$ be a nowhere
vanishing holomorphic $2$-form on $X_{I}$ such that
$\eta\wedge\bar{\eta}=2\gamma_{I}^{2}$.
Then there exist complex structures 
$J,K\in{\cal C}_{g}$ satisfying
\newline{$(1)$}
$(I,J,K)$ is a hyper-K\"ahler structure on $(X,g)$ 
with $\eta=\gamma_{J}+\sqrt{-1}\gamma_{K}$;
\newline{$(2)$}
${\cal H}^{2}_{+}(g)$ is a $3$-dimensional real vector 
space spanned by 
$\{\gamma_{I},\gamma_{J},\gamma_{K}\}$;
\newline{$(3)$}
${\cal C}_{g}=\{aI+bJ+cK;\,(a,b,c)\in{\Bbb R}^{3},\,
a^{2}+b^{2}+c^{2}=1\}$.
\end{thm}

\begin{pf}
See \cite[Sect.\,2, (i) $\Leftrightarrow$ (iii)]{H} 
for (1) and (2). 
Let $I'\in{\cal C}_{g}$. Since
$\gamma_{I'}\in{\cal H}^{2}_{+}(g)$ by 
\cite[Sect.\,2, (i) $\Leftrightarrow$ (iii)]{H}, 
we can write 
$\gamma_{I'}=a\gamma_{I}+b\gamma_{J}+c\gamma_{K}$
for some $a,b,c\in{\Bbb R}$. We get
$a^{2}+b^{2}+c^{2}=1$ by the relations
$\gamma_{I'}^{2}=\gamma_{I}^{2}=2dV_{g}$, 
$\gamma\wedge\eta=0$, and
$\eta\wedge\bar{\eta}=2\gamma_{I}^{2}$.
\end{pf}

\begin{lem}
Let $(I,J,K)$ be a hyper-K\"ahler structure on $(X,g)$.
The map from $SO(3)$ to the set of all hyper-K\"ahler
structures on $(X,g)$ defined by
$$
A=(a_{ij})\mapsto(a_{11}I+a_{12}J+a_{13}K,\,
a_{21}I+a_{22}J+a_{23}K,\,a_{31}I+a_{32}J+a_{33}K)
$$
is a bijection.
\end{lem}

\begin{pf}
It is obvious that the map defined as above 
is injective. Let $(I',J',K')$ be an arbitrary 
hyper-K\"ahler structure on $(X,g)$. 
By Theorem 2.3 (3), there is a real $3\times3$
matrix $B=(b_{ij})$ with
$$
I'=b_{11}I+b_{12}J+b_{13}K,
\quad
J'=b_{21}I+b_{22}J+b_{23}K,
\quad
K'=b_{31}I+b_{32}J+b_{33}K.
$$ 
We get $B\in SO(3)$ by the relations
$(I')^{2}=(J')^{2}=(K')^{2}=-1_{TX}$ and
$I'J'=-J'I'=K'$. This proves the surjectivity.
\end{pf}

By Lemma 2.4, the element 
$\gamma_{I}\wedge\gamma_{J}\wedge\gamma_{K}\in
\det{\cal H}^{2}_{+}(g)$ is independent 
of the choice of a hyper-K\"ahler structure 
$(I,J,K)$ on $(X,g)$, and
it defines an orientation on ${\cal H}^{2}_{+}(g)$.
In this note, ${\cal H}^{2}_{+}(g)$ is
equipped with this orientation.
\par
Let $A^{p}(X)$ denote the real vector space of 
real $C^{\infty}$ $p$-forms on $X$. 
For a complex structure $I$ on $X$, 
$A^{p,q}(X_{I})$ denotes the complex vector space of 
$C^{\infty}$ $(p,q)$-forms on $X_{I}$, and
$\Omega_{X_{I}}^{p}$ denotes the sheaf of holomorphic
$p$-forms on $X_{I}$.
\par
Recall that the $L^{2}$-inner product on 
$A^{p}(X)$ with respect to $g$ is defined by
$$
(f,f')_{L^{2}}
:=
\int_{X}f\wedge *_{g}f'
=
\int_{X}\langle f,f'\rangle_{x}\,dV_{g}(x),
\qquad
f,f'\in A^{p}(X).
$$
Equipped with the restriction of 
$(\cdot,\cdot)_{L^{2}}$,
${\cal H}^{2}_{+}(g)$ is a metrized vector space. 
Then
$\{\gamma_{I}/\sqrt{2},\gamma_{J}/\sqrt{2},
\gamma_{K}/\sqrt{2}\}$ is an oriented orthonormal 
basis of ${\cal H}^{2}_{+}(g)$
for every hyper-K\"ahler structure $(I,J,K)$ 
on $(X,g)$,
because $\gamma=\gamma_{I}\in A^{1,1}(X_{I})$ and 
$\eta=\gamma_{J}+\sqrt{-1}\gamma_{K}\in
H^{0}(X_{I},\Omega^{2}_{X_{I}})$ satisfy
the equations $\gamma\wedge\eta=\eta^{2}=0$.

\begin{lem}
The map from the set of hyper-K\"ahler structures on
$(X,g)$ to the set of oriented orthonormal basis of
${\cal H}^{2}_{+}(g)$ defined by
$(I,J,K)\mapsto
\left\{\gamma_{I}/\sqrt{2},\gamma_{J}/\sqrt{2},
\gamma_{K}/\sqrt{2}\right\}$, is a bijection.
\end{lem}

\begin{pf}
The result is an immediate consequence of Lemma 2.4.
\end{pf}


\section
{\bf Hyperbolic involutions on $K3$ surfaces
and Ricci-flat metrics}
\par
In this section, we recall a trick of Donaldson
that relates real $K3$ surfaces and $K3$ surfaces 
with anti-symplectic holomorphic involution.
We follow \cite[Chap.\,6, Sect.\,15]{DIK00}
and \cite[Sect.\,2 pp.21-22]{D2}.

\subsection
{}
{\bf Hyperbolic Involution}
\par
For a $C^{\infty}$ involution $\iota$ on $X$, we set
$$
H^{2}_{\pm}(X,{\Bbb Z})
:=
\{l\in H^{2}(X,{\Bbb Z});\,
\iota^{*}(l)=\pm l\},
\qquad
r(\iota)
:=
{\rm rank}_{\Bbb Z}\,H^{2}_{+}(X,{\Bbb Z}).
$$
By \cite[Cor.\,1.5.2]{Ni1}, 
$H^{2}_{+}(X,{\Bbb Z})\subset H^{2}(X,{\Bbb Z})$
is primitive and $2$-elementary.

\begin{defn}
A $C^{\infty}$ involution $\iota\colon X\to X$ is 
{\it hyperbolic} if the following two conditions 
are satisfied:
\newline{(1)}
$H^{2}_{+}(X,{\Bbb Z})$ has signature $(1,r(\iota)-1)$;
\newline{(2)}
$\iota$ is holomorphic with respect to 
a complex structure on $X$. 
\end{defn}

\begin{rem}
The second condition of Definition 3.1 does not 
seem very natural. We do no know if it is deduced 
from the first condition. Are there any 
$C^{\infty}$ involution on $X$ which is never 
holomorphic with respect to any complex structure 
on $X$, such that the invariant lattice 
$H^{2}_{+}(X,{\Bbb Z})$ is hyperbolic?
\end{rem}

\begin{defn}
For a hyperbolic involution $\iota\colon X\to X$, 
set
$$
{\cal E}^{\iota}
:=
\{g\in{\cal E};\,\iota^{*}g=g\}.
$$
\end{defn}

\begin{prop}
For every hyperbolic involution 
$\iota\colon X\to X$, one has
${\cal E}^{\iota}\not=\emptyset$. 
\end{prop}

\begin{pf}
There exists a complex structure $I$ on $X$
such that $\iota$ is holomorphic with respect to $I$.
Since $X_{I}$ is K\"ahler, there exists 
an $\iota$-invariant K\"ahler class $\kappa$ 
on $X_{I}$. Let $\gamma$ be the
unique Ricci-flat K\"ahler form representing $\kappa$.
Then $\iota^{*}\gamma=\gamma$ by the uniqueness 
of $\gamma$.
Let $g$ be the K\"ahler metric on $X$ 
whose K\"ahler form is $\gamma$. 
Then $g$ is Ricci-flat and $\iota$-invariant.
\end{pf}

Let $\iota\colon X\to X$ be a hyperbolic involution,
and let $g\in{\cal E}^{\iota}$. 
Then $\iota$ preserves ${\cal H}^{2}_{+}(g)$.
By identifying a real harmonic $2$-form on $(X,g)$ 
with its cohomology class in 
$H^{2}(X,{\Bbb R})$, we regard ${\cal H}^{2}_{+}(g)$ 
as an oriented subspace of $H^{2}(X,{\Bbb R})$. 
Since $*_{g}=1$ on ${\cal H}^{2}_{+}(g)$,
the cup-product $\langle\cdot,\cdot\rangle$ is 
positive-definite on 
${\cal H}^{2}_{+}(g)\subset H^{2}(X,{\Bbb R})$.

\begin{prop}
The orientation on ${\cal H}^{2}_{+}(g)$ is
preserved by $\iota$.
\end{prop}

\begin{pf}
Since $\iota$ is a diffeomorphism of $X$,
the result follows from \cite[Prop.\,6.2]{D1}.
\end{pf}

\begin{prop}
$(1)$ 
There exists a hyper-K\"ahler structure $(I,J,K)$ 
on $(X,g)$ with
\begin{equation}
\iota_{*}I=I\iota_{*},
\qquad
\iota_{*}J=-J\iota_{*},
\qquad
\iota_{*}K=-K\iota_{*}.
\end{equation}
$(2)$ 
If $(I',J',K')$ is another hyper-K\"ahler structure 
satisfying $(3.1)$, then there exists 
$\psi\in{\Bbb R}$ satisfying one of 
the following two equations:
\begin{equation}
(I',J',K')=
\begin{cases}
(I,\,\,\cos\psi J-\sin\psi K,\,\,
\sin\psi J+\cos\psi K),\\
(-I,\,\,\cos\psi J+\sin\psi K,\,\,
\sin\psi J-\cos\psi K).
\end{cases}
\end{equation}
\end{prop}

\begin{pf}
Set $\varPi(g)_{\pm}:=
\{\gamma\in{\cal H}^{2}_{+}(g);\,
\iota^{*}\gamma=\pm\gamma\}$.
Since the cup-product is positive definite on 
${\cal H}^{2}_{+}(g)$, the hyperbolicity of $\iota$
implies that $\dim\varPi(g)_{+}\leq1$. 
Since $\det\iota^{*}|_{{\cal H}^{2}_{+}(g)}=1$
by Proposition 3.5, 
we get $\dim\varPi(g)_{+}=1$ and $\dim\varPi(g)_{-}=2$.
Since $\iota$ is an involution preserving
the $L^{2}$-inner product $(\cdot,\cdot)_{L^{2}}$, 
$\iota^{*}$ is symmetric with respect to 
$(\cdot,\cdot)_{L^{2}}$. 
Hence there exists an oriented orthonormal basis 
$\{\gamma_{1},\gamma_{2},\gamma_{3}\}\subset
{\cal H}^{2}_{+}(g)$ with
\begin{equation}
\iota^{*}\gamma_{1}=\gamma_{1},
\qquad
\iota^{*}\gamma_{2}=-\gamma_{2},
\qquad
\iota^{*}\gamma_{3}=-\gamma_{3}.
\end{equation}
By Lemma 2.5, there exists a hyper-K\"ahler
structure $(I,J,K)$ on $(X,g)$ satisfying 
$\gamma_{1}=\gamma_{I}/\sqrt{2}$,
$\gamma_{2}=\gamma_{J}/\sqrt{2}$,
$\gamma_{3}=\gamma_{K}/\sqrt{2}$.
These equations, together with (2.1), (3.3) and
$\iota^{*}g=g$, yields (3.1). This proves (1). 
\par
Since $\dim\varPi(g)_{+}=1$, 
there exists $l\in{\Bbb R}$ such that 
$\gamma_{I'}=l\gamma_{I}$. 
This, together with
$\gamma_{I}^{2}=\gamma_{I'}^{2}=2dV_{g}$, 
implies that $I'=\pm I$. Since 
$\{\omega_{J}/\sqrt{2},\omega_{K}/\sqrt{2}\}$ and
$\{\omega_{J'}/\sqrt{2},\omega_{K'}/\sqrt{2}\}$ 
are orthonormal bases of $\varPi(g)_{-}$, 
there exists $\psi\in{\Bbb R}$ with
$$
(J',K')=(\cos\psi J\mp\sin\psi K,\,
\sin\psi J\pm\cos\psi K).
$$ 
Since $J'K'=I$ when $I'=I$ and 
since $J'K'=-I$ when $I'=-I$, 
we get (3.2).
\end{pf}

\begin{defn}
A hyper-K\"ahler structure $(I,J,K)$ on $(X,g)$ is
{\it compatible} with $\iota$ if Eq.\,(3.1) holds.
\end{defn}

\subsection
{\bf $2$-elementary $K3$ surfaces}
\par
Let $Y$ be a $K3$ surface, 
and let $\theta\colon Y\to Y$ be a holomorphic 
involution. 
Then $\theta$ is {\it anti-symplectic} if
\begin{equation}
\theta^{*}\eta=-\eta,
\qquad
\forall\,\eta\in H^{0}(Y,\Omega^{2}_{Y}).
\end{equation}

\begin{defn}
A $K3$ surface equipped with an anti-symplectic 
holomorphic involution is called 
a {\it $2$-elementary $K3$ surface}. 
\end{defn}

\begin{prop}
Let $(Y,\theta)$ be a $2$-elementary $K3$ surface 
equipped with a $\theta$-invariant Ricci-flat 
K\"ahler metric $g$.
Let $I$ be the complex structure on $X$ 
such that $Y=X_{I}$,
let $\eta$ be a holomorphic $2$-form on $Y$ 
such that $\eta\wedge\bar{\eta}=2\gamma_{I}^{2}$, 
and let $J,K\in{\cal C}_{g}$ be 
the complex structures such that
$\gamma_{J}={\rm Re}(\eta)$ and
$\gamma_{K}={\rm Im}(\eta)$.
Then 
\newline{$(1)$}
$\theta$ is a hyperbolic involution and
$g\in{\cal E}^{\theta}$;
\newline{$(2)$}
the hyper-K\"ahler structure $(I,J,K)$ on $(X,g)$ 
is compatible with $\theta$.
\end{prop}

\begin{pf}
By (3.4) and the $\theta$-invariance 
of $\gamma_{I}$, we get (3.1). 
The hyperbolicity of $\theta$ follows from 
e.g. \cite{DIK00}, \cite{Ni1}, 
\cite[Lemma 1.3 (1)]{Yo}.
\end{pf}

We refer to \cite{DIK00}, \cite{Ni2}, \cite{Yo} 
for more details about $2$-elementary $K3$ surfaces.

\subsection
{}
{\bf Real $K3$ surfaces}
\par
After \cite{DIK00}, \cite{Kh},
\cite[Sect.\,2 and Sect.\,3]{Ni3}, 
we make the following:

\begin{defn}
A $K3$ surface equipped with 
an {\it anti-holomorphic} involution 
is called a {\it real $K3$ surface}.
A point of a real $K3$ surface is {\it real}
if it is fixed by the anti-holomorphic involution.
\end{defn}

\begin{exmp}
Let $Y$ be an algebraic $K3$ surface defined 
over $\Bbb R$. 
Then there exists a projective embedding 
$j\colon Y\hookrightarrow{\Bbb P}^{N}({\Bbb C})$
defined over $\Bbb R$. The complex conjugation
${\Bbb P}^{N}({\Bbb C})\ni(z_{1}:\cdots:z_{N})\to
(\bar{z}_{1}:\cdots:\bar{z}_{N})\in
{\Bbb P}^{N}({\Bbb C})$ acts on $Y$ 
as an anti-holomorphic involution. 
Let $\sigma\colon Y\to Y$ be the involution 
induced by the complex conjugation 
on ${\Bbb P}^{N}({\Bbb C})$. 
Then the pair $(Y,\sigma)$ is a real $K3$ surface.
We refer to \cite{DIK00}, \cite{Kh},
\cite{Ni1}, \cite[Sect.\,2]{Ni3}
for more details about this example.
\end{exmp}

Let $(Y,\sigma)$ be a real $K3$ surface. 
Let $g$ be a K\"ahler
metric on $Y$ with K\"ahler form $\gamma$. 
Then $\sigma^{*}g$ is a K\"ahler metric 
with K\"ahler form $-\sigma^{*}\gamma$. 
Indeed, if $Y=X_{J}$, we get
\begin{equation}
(\sigma^{*}g)(J(u),v)
=
g(\sigma_{*}J(u),\sigma_{*}(v))
=
-g(J\sigma_{*}(u),\sigma_{*}(v))
=
-(\sigma^{*}\gamma)(u,v)
\end{equation}
for all $u,v\in TX$.
Hence $Y$ admits a $\sigma$-invariant K\"ahler metric
e.g. $g+\sigma^{*}g$.
By (3.5), the K\"ahler form and the K\"ahler class 
of a $\sigma$-invariant K\"ahler metric are 
anti-invariant with respect to the $\sigma$-action. 
In particular, there exists a K\"ahler class $\kappa$ 
on $Y$ with $\sigma^{*}\kappa=-\kappa$.

\begin{lem}
$(1)$ 
There exists  
$\eta\in H^{0}(Y,\Omega^{2}_{Y})\setminus\{0\}$ 
with
\begin{equation}
\sigma^{*}\eta=\bar{\eta}.
\end{equation}
$(2)$ 
Let $\kappa$ be a K\"ahler class on $Y$ with
$\sigma^{*}\kappa=-\kappa$, and let $\gamma$ be 
the Ricci-flat K\"ahler form representing $\kappa$. 
Then
\begin{equation}
\sigma^{*}\gamma=-\gamma.
\end{equation}
$(3)$ 
There exists a $\sigma$-invariant Ricci-flat 
K\"ahler metric on $Y$.
\end{lem}

\begin{pf}
{\bf(1)}
Let $\xi$ be a nowhere vanishing holomorphic 
$2$-form on $Y$. Since $\sigma$ is anti-holomorphic, 
$\sigma^{*}\bar{\xi}$ is a holomorphic $2$-form 
on $Y$. Then either $\xi+\sigma^{*}\bar{\xi}$ or
$(\xi-\sigma^{*}\bar{\xi})/\sqrt{-1}$ 
is a nowhere vanishing holomorphic $2$-form on $Y$ 
satisfying (3.6).
\newline{\bf (2)}
Let $g$ be the Riemannian metric on $Y$ whose K\"ahler
form is $\gamma$. By (3.5), $-\sigma^{*}\gamma$ 
is the K\"ahler form of $\sigma^{*}g$ 
representing $\kappa$. 
By the Ricci-flatness of $\gamma$,
there exists a real non-zero constant $C$ with
$C\,\gamma^{2}=\eta\wedge\bar{\eta}$.
This, together with (3.6), yields that
$$
C\,(-\sigma^{*}\gamma)^{2}
=
\sigma^{*}\eta\wedge\sigma^{*}\bar{\eta}
=
\bar{\eta}\wedge\eta=\eta\wedge\bar{\eta}.
$$
This implies the Ricci-flatness of 
$-\sigma^{*}\gamma$.
By the uniqueness of the Ricci-flat K\"ahler form
in the K\"ahler class $\kappa$, we get (3.7).
\newline{\bf (3)}
By (2), there exists a Ricci-flat K\"ahler metric 
$g$ on $Y$ whose K\"ahler form satisfies (3.7). 
Since $\sigma$ is anti-holomorphic, 
we get $\sigma^{*}g=g$ by (3.7).
\end{pf}

\begin{defn}
A holomorphic $2$-form $\eta$ on 
a real $K3$ surface $(Y,\sigma)$ is  
{\it defined over $\Bbb R$} if Eq.\,$(3.6)$ holds.
\end{defn}

\begin{prop}
Let $(Y,\sigma)$ be a real $K3$ surface 
equipped with a $\sigma$-invariant Ricci-flat 
K\"ahler metric $g$.
Let $J$ be the complex structure on $X$ 
with $Y=X_{J}$,
let $\eta$ be a holomorphic $2$-form on $Y$ 
defined over $\Bbb R$ with 
$\eta\wedge\bar{\eta}=2\gamma_{J}^{2}$, 
and let $I,K\in{\cal C}_{g}$ 
be the complex structures with
$\gamma_{I}=-{\rm Re}\,\eta$ and
$\gamma_{K}={\rm Im}\,\eta$.
Then
\newline{$(1)$}
$\sigma$ is a hyperbolic involution and 
$g\in{\cal E}^{\sigma}$;
\newline{$(2)$}
the hyper-K\"ahler structure $(I,J,K)$ is
compatible with $(g,\sigma)$.
\end{prop}

\begin{pf}
By (3.6) and (3.7), we get
\begin{equation}
\sigma^{*}\gamma_{I}=\gamma_{I},
\qquad
\sigma^{*}\gamma_{J}=-\gamma_{J},
\qquad
\sigma^{*}\gamma_{K}=-\gamma_{K},
\end{equation}
which, together with $\sigma^{*}g=g$, implies (3.1).
Hence it suffices to verify the hyperbolicity 
of $\sigma$.
Consider the $K3$ surface $X_{I}$. 
By (3.1) and (3.8), $\sigma\colon X_{I}\to X_{I}$ 
is an anti-symplectic holomorphic involution. 
Hence $\sigma$ is hyperbolic.
\end{pf}

\begin{prop}
Let $\iota\colon X\to X$ be a hyperbolic involution, 
and let $g\in{\cal E}^{\iota}$. Let $(I,J,K)$ be 
a hyper-K\"ahler structure on $(X,g)$ 
compatible with $\iota$. Then
\newline{$(1)$}
$(X_{I},\iota)$ is a $2$-elementary $K3$ surface, 
and $\gamma_{J}+\sqrt{-1}\gamma_{K}$ is 
a holomorphic $2$-form on $X_{I}$;
\newline{$(2)$}
$(X_{J},\iota)$ is a real $K3$ surface, and
$\gamma_{I}+\sqrt{-1}\gamma_{K}$ is a holomorphic
$2$-form on $X_{J}$ defined over $\Bbb R$.
\end{prop}

\begin{pf}
The result follows from (3.1) and Propositions 3.9 
and 3.14.
\end{pf}

\subsection
{}
{\bf The period map for Ricci-flat metrics compatible
with involution}
\par
Let $M\subset{\Bbb L}_{K3}$ be a sublattice.

\begin{defn}
A hyperbolic involution $\iota\colon X\to X$ is 
{\it of type $M$} if there exists an isometry of
lattices 
$\alpha\colon H^{2}(X,{\Bbb Z})\cong{\Bbb L}_{K3}$
such that $M=\alpha(H^{2}_{+}(X,{\Bbb Z}))$. 
An isometry 
$\alpha\colon H^{2}(X,{\Bbb Z})\cong{\Bbb L}_{K3}$
with this property is called 
a {\it marking of type $M$}.
\end{defn}

Let $\iota$ be a hyperbolic involution of type $M$, 
and let 
$\alpha\colon H^{2}(X,{\Bbb Z})\cong{\Bbb L}_{K3}$ 
be a marking of type $M$. Then
$M\subset{\Bbb L}_{K3}$ is a primitive, 
$2$-elementary, hyperbolic sublattice by 
\cite[Cor\,1.5.2]{Ni1}.
The orthogonal complement of $M$ in ${\Bbb L}_{K3}$
is denoted by $M^{\perp}$. Then 
$M^{\perp}=\alpha(H^{2}_{-}(X,{\Bbb Z}))$.
We set $r(M):={\rm rank}_{\Bbb Z}\,M$ and
$$
\Omega_{M}
:=
\{[\eta]\in{\Bbb P}(M^{\perp}\otimes{\Bbb C});\,
\langle\eta,\eta\rangle=0,\,
\langle\eta,\bar{\eta}\rangle>0\}.
$$ 
Since $M^{\perp}$ has signature $(2,20-r(M))$, 
$\Omega_{M}$ consists of two connected components,
each of which is isomorphic to a symmetric bounded 
domains of type IV of dimension $20-r(M)$ 
(cf. \cite[p.282, Lemma 20.1]{BPV}). 
Then $\Omega_{M}$ is the period domain 
for $2$-elementary $K3$ surfaces of type $M$ 
by \cite[Sect.\,1.4]{Yo}. 
Notice that the two connected components 
of $\Omega_{M}$ are exchanged by
the complex conjugation on 
${\Bbb P}(M^{\perp}\otimes{\Bbb C})$.

\begin{lem}
Let $\iota\colon X\to X$ be a hyperbolic involution 
of type $M$, and let $\alpha$ be a marking 
of type $M$.
Let $g\in{\cal E}^{\iota}$, and 
let $(I,J,K)$ be a hyper-K\"ahler structure 
on $(X,g)$ compatible with $\iota$. 
Then the pair of conjugate points
$[\alpha(\gamma_{J}\pm\sqrt{-1}\gamma_{K})]
\in\Omega_{M}$ is independent of 
the choice of $(I,J,K)$ compatible with $\iota$.
\end{lem}

\begin{pf}
By Proposition 3.15 (1),  
$[\alpha(\gamma_{J}+\sqrt{-1}\gamma_{K})]$ is 
the period of a marked $2$-elementary $K3$ surface 
of type $M$. Hence 
$[\alpha(\gamma_{J}+\sqrt{-1}\gamma_{K})]
\in\Omega_{M}$ by \cite[Sect.\,1.4]{Yo}.
Since the complex conjugation preserves $\Omega_{M}$, 
we get
$[\alpha(\gamma_{J}\pm\sqrt{-1}\gamma_{K})]
\in\Omega_{M}$.
\par
Let $(I',J',K')$ be an arbitrary hyper-K\"ahler 
structure on $(X,g)$ compatible with $\iota$.
By Proposition 3.6 (2), there exists 
$\psi\in{\Bbb R}$ such that
$$
\gamma_{J'}+\sqrt{-1}\gamma_{K'}
=
e^{\sqrt{-1}\psi}(\gamma_{J}\pm\sqrt{-1}\gamma_{K}).
$$
Hence $[\alpha(\gamma_{J}\pm\sqrt{-1}\gamma_{K})]=
[\alpha(\gamma_{J'}\pm\sqrt{-1}\gamma_{K'})]
\in\Omega_{M}$.
\end{pf}

\begin{defn}
With the same notation as in Lemma 3.17,
the pair of conjugate points
$[\alpha(\gamma_{J}\pm\sqrt{-1}\gamma_{K})]
\in\Omega_{M}$ is called the {\it period} 
of $(g,\alpha)$ and is denoted by
$$
\varpi_{M}(g,\alpha):=
[\alpha(\gamma_{J}\pm\sqrt{-1}\gamma_{K})].
$$
\end{defn}


\section
{\bf An invariant of Ricci-flat metric compatible
with involution}
\par
Throughout this section, 
we fix the following notation.
Let $\iota\colon X\to X$ be a hyperbolic involution 
of type $M$, and let 
$\alpha\colon H^{2}(X,{\Bbb Z})\cong{\Bbb L}_{K3}$ 
be a marking of type $M$.
Let ${\Bbb Z}_{2}=\langle\iota\rangle$ be the group 
of diffeomorphisms of $X$ generated by $\iota$. 
Let $g\in{\cal E}^{\iota}$.

\subsection
{}
{\bf Equivariant determinant of the Laplacian}
\par
Let $d^{*}\colon A^{1}(X)\to C^{\infty}(X)$ 
be the formal adjoint of the exterior
derivative $d\colon C^{\infty}(X)\to A^{1}(X)$ 
with respect to the $L^{2}$-inner product 
induced by $g$.
The Laplacian of $(X,g)$ is defined as
$\Delta_{g}=\frac{1}{2}d^{*}d$. We define
$$
C^{\infty}_{\pm}(X):=
\{f\in C^{\infty}(X);\,\iota^{*}f=\pm f\}.
$$
Since $\iota$ preserves $g$, $\Delta_{g}$ commutes 
with the $\iota$-action on $C^{\infty}(X)$. 
Hence $\Delta_{g}$ preserves 
the subspaces $C^{\infty}_{\pm}(X)$. 
We set 
$$
\Delta_{g,\pm}:=\Delta_{g}|_{C^{\infty}_{\pm}(X)}.
$$
Define the spectral zeta function
of $\Delta_{g,\pm}$ as
$$
\zeta_{g,\pm}(s)
:=
{\rm Tr}\left\{\Delta_{g,\pm}
|_{(\ker\Delta_{g})^{\perp}}\right\}^{-s}
=
{\rm Tr}\left[\frac{1\pm\iota^{*}}{2}\circ
\left(\Delta_{g}
|_{(\ker\Delta_{g})^{\perp}}\right)^{-s}\right],
\qquad
{\rm Re}\,s\gg0.
$$
Then $\zeta_{g,\pm}(s)$ converges absolutely
for ${\rm Re}\,s\gg0$, it extends meromorphically to
the complex plane $\Bbb C$, and it is holomorphic 
at $s=0$.

\begin{defn}
(1) The equivariant determinant of $\Delta_{g}$ with
respect to ${\Bbb Z}_{2}=\langle\iota\rangle$ is 
defined by
$$
{\rm det}_{{\Bbb Z}_{2}}^{*}\Delta_{g}(\iota)
:=
\exp[-\zeta'_{g,+}(0)+\zeta'_{g,-}(0)].
$$
\newline{(2)}
For a real $K3$ surface $(Y,\sigma)$ and 
a $\sigma$-invariant Ricci-flat K\"ahler metric $g$, 
set
$$
{\rm det}_{{\Bbb Z}_{2}}^{*}\Delta_{Y,g}(\sigma)
:=
{\rm det}_{{\Bbb Z}_{2}}^{*}\Delta_{g}(\sigma).
$$
\end{defn}

\subsection
{\bf Equivariant determinant of the Laplacian
and equivariant analytic torsion}
\par

Let $(I,J,K)$ be a hyper-K\"ahler structure 
on $(X,g)$ compatible with $\iota$. 
By Proposition 3.15 (1), $\iota$ is 
an anti-symplectic holomorphic involution 
on $X_{I}$.
\par
Let $\square_{g,I}^{0,q}$ be 
the $\bar{\partial}$-Laplacian
acting on $(0,q)$-forms on the K\"ahler manifold 
$(X_{I},\gamma_{I})$. 
By the definition of $\Delta_{g}$
and the K\"ahler identities, one has
$\Delta_{g}=\square_{g,I}^{0,0}$. We set
$$
\zeta^{0,q}(g,I,\iota)(s)
:=
{\rm Tr}\left[\iota^{*}(\square_{g,I}^{0,q}
|_{(\ker\square_{g,I}^{0,q})^{\perp}})^{-s}\right],
\qquad
{\rm Re}\,s\gg0.
$$
Then 
\begin{equation}
\zeta^{0,1}(g,I,\iota)(s)
=
\zeta^{0,0}(g,I,\iota)(s)+
\zeta^{0,2}(g,I,\iota)(s),
\end{equation}
\begin{equation}
\zeta^{0,0}(g,I,\iota)(s)
=
\zeta^{+}_{g}(s)-\zeta^{-}_{g}(s).
\end{equation}
After \cite{Bi} and \cite{KR}, 
we make the following:

\begin{defn}
The equivariant analytic torsion of
$(X_{I},\gamma_{I},\iota)$ is defined by
$$
\tau_{{\Bbb Z}_{2}}(g,I,\iota)
:=
\exp\left[\zeta^{0,1}(g,I,\iota)'(0)-
2\zeta^{0,2}(g,I,\iota)'(0)\right].
$$
\end{defn}

\begin{lem}
The following identity holds
$$
\tau_{{\Bbb Z}_{2}}(g,I,\iota)=
\left({\rm det}_{{\Bbb Z}_{2}}^{*}
\Delta_{g}(\iota)\right)^{-2}.
$$
\end{lem}

\begin{pf}
Let $K_{X_{I}}$ be the canonical line bundle 
of $X_{I}$, and set 
$\eta_{I}=\gamma_{J}+\sqrt{-1}\gamma_{K}\in
H^{0}(X_{I},K_{X_{I}})$. Since $\gamma_{J}$ and 
$\gamma_{K}$ are parallel with respect
to the Levi-Civita connection of $(X,g)$, 
so is $\eta_{I}$. 
The isomorphism of complex line bundles 
$\otimes\bar{\eta}\colon{\cal O}_{X_{I}}\cong
\overline{K}_{X_{I}}$ induces an isometry 
with respect to the $L^{2}$-inner products:
$$
\otimes\bar{\eta}/\sqrt{2}\colon
C^{\infty}(X)\ni f\to 
f\cdot\bar{\eta}/\sqrt{2}\in A^{0,2}(X_{I}).
$$ 
\par
Let $E_{g}(\lambda)$ 
(resp. $E_{g,I}^{0,2}(\lambda)$)
be the eigenspace of $\Delta_{g}$ 
(resp. $\square_{g,I}^{0,2}$)
with respect to the eigenvalue
$\lambda\in{\Bbb R}$. Then $\iota$ preserves
$E_{g}(\lambda)$ and $E_{g,I}^{0,2}(\lambda)$.
Let $E_{g}(\lambda)_{\pm}$ and 
$E_{g,I}^{0,2}(\lambda)_{\pm}$ be 
the $\pm1$-eigenspaces of the $\iota$-actions on
$E_{g}(\lambda)$ and $E_{g,I}^{0,2}(\lambda)$,
respectively.
Since $\iota^{*}\bar{\eta}=-\bar{\eta}$ and
$$
\square^{0,2}_{g,I}(f\cdot\bar{\eta})=
(\Delta_{g} f)\cdot\bar{\eta},
\qquad
f\in C^{\infty}(X),
$$
we get the isomorphism
$\otimes\bar{\eta}/\sqrt{2}\colon
E_{g}(\lambda)_{\pm}\cong 
E^{0,2}_{g,I}(\lambda)_{\mp}$ 
for all $\lambda\in{\Bbb R}$, which yields that 
\begin{equation}
\zeta^{0,2}(g,I,\iota)(s)=
-\zeta^{+}_{g}(s)+\zeta^{-}_{g}(s),
\qquad
s\in{\Bbb C}.
\end{equation}
By (4.1), (4.2) and (4.3), we get
\begin{equation}
\begin{aligned}
\log\tau_{{\Bbb Z}_{2}}(g,I,\iota)
&
=
\zeta^{0,1}(g,I,\iota)'(0)-
2\zeta^{0,2}(g,I,\iota)'(0)
\\
&
=
\zeta^{0,0}(g,I,\iota)'(0)-
\zeta^{0,2}(g,I,\iota)'(0)
\\
&
=
2\left.\frac{d}{ds}\right|_{s=0}
(\zeta^{+}_{g}(s)-\zeta^{-}_{g}(s))
=
-2\log{\rm det}_{{\Bbb Z}_{2}}^{*}\Delta_{g}(\iota).
\end{aligned}
\end{equation}
This completes the proof of Lemma 4.3.
\end{pf}

\subsection
{}
{\bf A function $\tau_{\iota}$ 
on ${\cal E}^{\iota}$}
\par

Let $X^{\iota}$ be the set of fixed points of
$\iota$: 
$$
X^{\iota}:=\{x\in X;\,\iota(x)=x\}.
$$
By \cite[Th.\,3.10.6]{Ni1} or \cite[Th.\,4.2.2]{Ni2}, 
$X^{\iota}$ is either the empty set or
the disjoint union of finitely many compact, 
connected, orientable two-dimensional manifolds. 
Moreover, $r(\iota)=10$ when $X^{\iota}=\emptyset$.
\par
When $X^{\iota}\not=\emptyset$, the Riemannian metric
$g|_{X^{\iota}}$ induces a complex structure 
on $X^{\iota}$ such that 
$g|_{X^{\iota}}$ is K\"ahler.
Equipped with this complex structure, 
$X^{\iota}$ is a complex submanifold of $X_{I}$, 
since $\iota$ is holomorphic with respect to $I$.
Let
$$
X^{\iota}=\amalg_{i}C_{i}
$$
be the decomposition into the connected components.
Let $\Delta_{(C_{i},g|_{C_{i}})}:=\frac{1}{2}d^{*}d$ 
be the Laplacain of the Riemannian manifold
$(C_{i},g|_{C_{i}})$, and let 
$$
\zeta_{(C_{i},g|_{C_{i}})}(s)
:=
{\rm Tr}\left[\Delta_{(C_{i},g|_{C_{i}})}
|_{(\ker\Delta_{(C_{i},g|_{C_{i}})})^{\perp}}
\right]^{-s}
$$ 
be the spectral zeta function of 
$\Delta_{(C_{i},g|_{C_{i}})}$.
The regularized determinant of 
$\Delta_{(C_{i},g|_{C_{i}})}$ is defined as 
$$
{\rm det}^{*}\Delta_{(C_{i},g|_{C_{i}})}
:=
\exp\left(-\zeta_{(C_{i},g|_{C_{i}})}'(0)\right).
$$
\par
Similarly, let $\tau(C_{i,I},\gamma_{I}|_{C_{i}})$ 
be the analytic torsion of the trivial Hermitian 
line bundle on the K\"ahler manifold 
$(C_{i},I,\gamma_{I}|_{C_{i}})$ (cf. \cite{RS}).
For all $i$, one has
\begin{equation}
\tau(C_{i,I},\gamma_{I}|_{C_{i}})
=
({\rm det}^{*}\Delta_{(C_{i},g|_{C_{i}})})^{-1}.
\end{equation}
\par
We define a function
$\tau_{\iota}$ on ${\cal E}^{\iota}$ and a function
$\tau_{M}$ on the moduli space of $2$-elementary 
$K3$ surfaces of type $M$ (cf. \cite[Def.\,5.1]{Yo}) 
as follows:

\begin{defn}
Let $(I,J,K)$ be a hyper-K\"ahler structure 
on $(X,g)$ compatible with $\iota$.
When $X^{\iota}\not=\emptyset$, set
$$
\tau_{\iota}(g)
:=
\left\{{\rm det}_{{\Bbb Z}_{2}}^{*}
\Delta_{g}(\iota)\right\}^{-2}
\prod_{i}{\rm Vol}(C_{i},g|_{C_{i}})\,
({\rm det}^{*}\Delta_{(C_{i},g|_{C_{i}})})^{-1},
$$
$$
\tau_{M}(X_{I},\iota)
:=
\tau_{{\Bbb Z}_{2}}(X_{I},\gamma_{I})(\iota)
\prod_{i}{\rm Vol}(C_{i},\gamma_{I}|_{C_{i}})\,
\tau(C_{i,I},\gamma_{I}|_{C_{i}}).
$$
When $X^{\iota}=\emptyset$, set
$$
\tau_{\iota}(g)
:=
\left\{{\rm det}_{{\Bbb Z}_{2}}^{*}
\Delta_{g}(\iota)\right\}^{-2},
\qquad
\tau_{M}(X_{I},\iota)
:=
\tau_{{\Bbb Z}_{2}}(X_{I},\gamma_{I})(\iota).
$$
\end{defn}

Notice that $(X,g)$ has volume $1$ for
$g\in{\cal E}^{\iota}$.
By \cite[Th.\,5.7]{Yo}, $\tau_{M}(X_{I},\iota)$ is
independent of the choice of an $\iota$-invariant
Ricci-flat K\"ahler metric on $X_{I}$.

\begin{lem}
If the hyper-K\"ahler structure $(I,J,K)$ 
on $(X,g)$ is compatible with $\iota$, then
\begin{equation}
\tau_{\iota}(g)=\tau_{M}(X_{I},\iota).
\end{equation}
In particular, one has
\begin{equation}
\tau_{M}(X_{I},\iota)=\tau_{M}(X_{-I},\iota).
\end{equation}
\end{lem}

\begin{pf}
The first result follows from Lemma 4.3 and (4.5).
If $(I,J,K)$ is compatible with $\iota$, so is
$(-I,J,-K)$. Hence the second result follows from
the first one.
\end{pf}

In the next theorem, we shall use the notion of
automorphic forms on $\Omega_{M}$, for which we refer 
to \cite[Sect.\,3]{Yo}. 
For an automorphic form $\Psi$ on $\Omega_{M}$, 
its norm $\|\Psi\|$ is a function on $\Omega_{M}$ 
defined in \cite[Def.\,3.16]{Yo}. 
If $X^{\iota}=\emptyset$ or
if every connected component of $X^{\iota}$ 
is diffeomorphic to a $2$-sphere, 
then $\Psi$ is an automorphic form 
in the classical sense and $\|\Psi\|$ coincides with 
the Petersson norm of $\Psi$.

\begin{thm} 
There exist $\nu(M)\in{\Bbb N}$ and
an automorphic form $\Phi_{M}$ on $\Omega_{M}$ 
of weight $((r(M)-6)\nu(M),4\nu(M))$ 
for some cofinite subgroup of $O(M^{\perp})$ 
satisfying
\newline{$(1)$}
$\|\Phi_{M}([\eta])\|=\|\Phi_{M}([\overline{\eta}])\|$
for all $[\eta]\in\Omega_{M}$;
\newline{$(2)$}
For all $g\in{\cal E}^{\iota}$,
\begin{equation}
\tau_{\iota}(g)=
\|\Phi_{M}(\varpi_{M}(g,\alpha))
\|^{-\frac{1}{2\nu(M)}}.
\end{equation}
\end{thm}

\begin{pf}
Let $\Phi_{M}$ be the automorphic form 
as in \cite[Th.\,5.2]{Yo}.
Let $(I,J,K)$ be a hyper-K\"ahler structure 
on $(X,g)$ compatible with $\iota$.
Let $(X_{I},\iota)$ be a $2$-elementary $K3$ surface
of type $M$. Then so is $(X_{-I},\iota)$. 
Since an anti-holomorphic $2$-form on $X_{I}$ is 
a holomorphic $2$-form on $X_{-I}$,
the Griffiths period of $(X_{-I},\iota)$ 
in the sense of \cite[(1.11)]{Yo} is the complex
conjugate of the Griffiths period of $(X_{I},\iota)$.
This, together with \cite[Th.\,5.2]{Yo} and (4.7),
implies the first assertion.
Since $\varpi_{M}(g,\alpha)=
\alpha(\gamma_{J}\pm\sqrt{-1}\gamma_{K})$ and
since $\gamma_{J}+\sqrt{-1}\gamma_{K}\in
H^{0}(X_{I},\Omega^{2}_{X_{I}})$, the second 
assertion follows from \cite[Th.\,5.2]{Yo} 
and (4.6).
\end{pf}

We assume that $\iota$ has no fixed points. 
By Proposition 3.15 (1), $\iota$ is a holomorphic
involution on $X_{I}$ without fixed points,
so that the quotient
$X_{I}/\iota$ is an Enriques surface by
\cite[Chap.\,8, Lemma 15.1]{BPV}. 
By \cite[Chap.\,8, Lemma 19.1]{BPV}, 
there exists an isometry 
$\alpha\colon H^{2}(X,{\Bbb Z})\cong{\Bbb L}_{K3}$ 
such that
$$
\alpha\,\iota^{*}\alpha^{-1}(a,b,c,x,y)=(b,a,-c,y,x),
\qquad
a,b,c\in{\Bbb U},
\quad
x,y\in{\Bbb E}_{8}.
$$
Set ${\frak L}:=\alpha(H^{2}_{+}(X,{\Bbb Z}))$. 
Then $\iota$ is of type $\frak L$.
We refer to \cite[Chap.\,8, Sects.\,15-21]{BPV}
for more details about Enriques surfaces.
\par
Let $\Phi$ be the {\it Borcherds $\Phi$-function},
which is an automorphic form of weight $4$ 
on the period domain for Enriques surfaces 
by \cite{Bo1}. 
By \cite[Th.\,8.2]{Yo},
there exists a constant $C_{\frak L}\not=0$ 
such that 
\begin{equation}
\Phi_{\frak L}=C_{\frak L}\,\Phi.
\end{equation}
Since $\iota$ has no fixed points, 
we may choose $\nu({\frak L})=1$ in Theorem 4.6 
by the definition of $\nu(M)$ in \cite[pp.\,79]{Yo}.

\begin{cor} 
Let $(Y,\sigma)$ be a real $K3$ surface 
without real points.
Let $g$ be a $\sigma$-invariant Ricci-flat K\"ahler 
metric on $Y$ with volume $1$. Let $\omega_{g}$ be
the K\"ahler form of $g$, and 
let $\eta_{g}$ be a holomorphic $2$-form on $Y$ 
defined over $\Bbb R$ such that 
$\eta_{g}\wedge\bar{\eta}_{g}=2\omega_{g}^{2}$. 
Let $\alpha$ be a marking of type $\frak L$.
Under the identifications of $\omega_{g}$ and 
$\eta_{g}$ with their cohomology classes, 
the following identity holds:
$$
{\rm det}_{{\Bbb Z}_{2}}^{*}\Delta_{Y,g}(\sigma)
=
C_{\frak L}^{\frac{1}{4}}\,
\|\Phi([\alpha(\gamma_{g}+
\sqrt{-1}{\rm Im}\,\eta_{g})])\|^{\frac{1}{4}}.
$$
\end{cor}

\begin{pf}
By Proposition 3.14 and Definition 3.18, we get
$\varpi_{\frak L}(g,\alpha)=
[\alpha(\gamma_{g}+\sqrt{-1}{\rm Im}\,\eta_{g})]$.
Substituting this equality and (4.9) into (4.8),
we get the result.
\end{pf}


\end{document}